\documentclass[12pt]{article}

\usepackage[utf8]{inputenc}

\usepackage{amsfonts}
\usepackage{amssymb}
\usepackage{amsmath}
\usepackage{amsthm}
\usepackage[T2A]{fontenc}

\usepackage{amsfonts}
\usepackage{amssymb}
\usepackage{amsmath}
\usepackage{amsthm}

\begin{document} 

\begin{Large}
\centerline{A note on Diophatine approximation in $\rm{ SL}_2(\mathbb{R})$.}
\vskip+0.3cm
\centerline{by Nikolay Moshchevitin}
\end{Large}

 \vskip+1.0cm

 \vskip+1.0cm

{\bf 1. Two-dimensional lattices and $2\times 2$ matrices.}

Suppose that $\rho >0$.  For reals $A_1, B_1$ we define 
\begin{equation}\label{rho}
 A_2 = -\rho A_1, B_2 = -\rho B_1.
\end{equation}
Let $r \in (0,1)$ be a positive constant such that for any $A_1,B_1\asymp Q$  
the equation  
\begin{equation}\label{dee}
xy-zw =1
\end{equation}
has a solution $x,y,z,w\in \mathbb{Z}$ with 
\begin{equation}\label{ii}
A_1 -Q^r\le x\le A_1+Q^r,\,
B_1 -Q^r\le z \le B_1+Q^r,\,\,
A_2 -Q^r\le y\le A_2+Q^r,\,
B_2-Q^r\le w \le B_2+Q^r.
\end{equation}
Of course $r$ does not depend on $\rho$.
(We can take as $r$ any  number greater than
 $ 3/4$, see Section 5).

{\bf Theorem 1.}\,\,{\it
Suppose that we have two unimodular matrices
\begin{equation}\label{maa}
\left(
\begin{array}{cc}
\alpha &\gamma\cr
\beta&\delta
\end{array}\right)
 ,
\left(
\begin{array}{cc}
\xi_1 &\xi_2\cr
\eta_1&\eta_2
\end{array}\right) \in {\rm SL}_2(\mathbb{R}).
\end{equation}
Suppose that $\delta \neq 0$ and
$\beta/\delta \not\in \mathbb{Q}$.
 Then there exists a sequence of reals $t_\nu \to \infty, \, \nu \to\infty$
and a sequence of  unimodular integer matrices
$$
\left(
\begin{array}{cc}
l_1^1 &l_1^2\cr
l_2^1&l_2^2
\end{array}\right) =
\left(
\begin{array}{cc}
l_1^1 (\nu) &l_1^2(\nu)\cr
l_2^1(\nu)&l_2^2(\nu)
\end{array}\right) \in {\rm SL}_2 (\mathbb{Z})
$$
 such that
\begin{equation}\label{erre}
\left(
\begin{array}{cc}
1 &t_\nu\cr
0&1
\end{array}\right)
\left(
\begin{array}{cc}
\alpha &\gamma\cr
\beta&\delta
\end{array}\right)
\left(
\begin{array}{cc}
l_1^1 &l_1^2\cr
l_2^1&l_2^2
\end{array}\right)
-
\left(
\begin{array}{cc}
\xi_1 &\xi_2\cr
\eta_1&\eta_2
\end{array}\right)
= O(|t_\nu|^{\frac{r-1}{r+1}}),
\end{equation}
where the constant in $O(\cdot )$ may depend on the size of matrices (\ref{maa}).

}

{\bf Corollary.}\,\, {\it
Let $\underline{\Lambda}, \Gamma$ be two unimodular lattices in $\mathbb{R}^2$.  
Suppose that
$$
\underline{
\Lambda} =
\left(
\begin{array}{cc}
\alpha &\beta\cr
\gamma&\delta
\end{array}\right)
\mathbb{Z}^2.
$$
and
$\beta/\delta \not\in \mathbb{Q}$.
Then there exists a sequence of reals $t_\nu \to \infty, \, \nu \to\infty$ such that for the
lattice $\underline{\Lambda}_{\,t_\nu} =  
\underline{\Lambda} \cdot \left(
\begin{array}{cc}
1 &0\cr
t_\nu&1
\end{array}\right) 
$
from the orbit of the lattice $\underline{\Lambda}$
 in the horocycle flow on ${\rm SL}_2 (\mathbb{Z}) \setminus  {\rm SL}_2 (\mathbb{R})$  one has
$$
dist( \underline{\Lambda}_{\, t_\nu}, \Gamma ) = O(|t_\nu|^{\frac{r-1}{r+1}})
$$
(here we consider the natural distance in the space of lattices (see \cite{cas})).
}

{\bf Theorem 2.}\,\,{\it
Consider a function $\psi (t)$ decreasing to zero as $t\to +\infty$ and such that
$ \psi (t) = O(t^{-1})$.
Define $\rho (t) $ to be the function inverse to the function $t \mapsto 1/\psi (t)$.
Suppose that
under the conditions of Theorem 1 one has
\begin{equation}\label{zabo}
\min_{p\in \mathbb{Z}}
|(\beta /\delta)\cdot q-p| \ge \psi (q),\,\,\,\, \forall q = 1,2,3,.... .
\end{equation}
Then there exists a positive constant $C$ such that for any $T\ge 1$  there exists a solution 
$$
\left(
\begin{array}{cc}
l_1^1 &l_1^2\cr
l_2^1&l_2^2
\end{array}\right)   \in {\rm SL}_2 (\mathbb{Z})
,\,\,\,
t\in \mathbb{R}$$
 of the system
\begin{equation}\label{oop}
\left|\left|
\left(
\begin{array}{cc}
1 &t\cr
0&1
\end{array}\right)
\left(
\begin{array}{cc}
\alpha &\gamma\cr
\beta&\delta
\end{array}\right)
\left(
\begin{array}{cc}
l_1^1 &l_1^2\cr
l_2^1&l_2^2
\end{array}\right)
-
\left(
\begin{array}{cc}
\xi_1 &\xi_2\cr
\eta_1&\eta_2
\end{array}\right)
\right|\right|\le C T^{\frac{r}{1+r} }/\rho \left( T^{\frac{1}{1+r}}\right)
,\,\,\,\,\,
\,\,\,
1\le |t |\le T
.
\end{equation}
Here $||\cdot||$ stands for the maximum of absolute values of elements of a matrix.}

{\bf Example.} \, If $\psi (t) = t^{-\omega}, \, \omega \ge 1$ then $\rho(t)
 = t^{1/\omega}$ and the right hand side of (\ref{oop}) is equal to
$$
O\left(T^{\frac{r}{1+r}-\frac{1}{\omega(1+r)}}\right).
$$
This bound is non-trivial provided $ \omega < r^{-1}$.
Moreover in the case when $\beta/\delta$ is a badly approximable number (that is, $\psi (t) = \kappa t^{-1}$
with a positive $\kappa$) the
right hand side of (\ref{oop})  is equal to
$ T^{\frac{r-1}{r+1}}$; in this case we see that the result of Theorem 1 holds uniformly.

Note that when $r$
is close to $3/4$ then the  exponent $\frac{r-1}{r+1}$ is close to $-1/7$.

In the next section we discuss some important history.
We give all the proofs in Sections 3--5.

{\bf 2. Some history.}

Results similar to our Theorem 1 are known for a long time.
They are related to  quantitative version of the famous Ratner's orbit closure theorem, in the simplest case of unipotent flow on $\rm{ SL}_2(\mathbb{R})$.
We refer to a wonderful book \cite{baak}  and the bibliography therein, as well as to the oridinal paper
\cite{R1}
by M. Ratner.
The qualitive result of such a type for the simplest case of $\rm{ SL}_2(\mathbb{R})$
was known much earlier.

Certain results similar to our Theorem 1
are due to A. Str\"ombergsson \cite{SSS} 
and F. Maucourant and
B. Weiss \cite{MW}.
From their theorems  a result of the same form as our Theorem 1 follows immediately,
but in the right hand side of (\ref{erre}) then we have
$O(t_\nu^{-\delta})$ with an effective positive $\delta$ which is not calculated explicitly. We do not compare our exponent to those from \cite{SSS,MW}.

However papers  \cite{SSS,MW} rely on the methods of dynamical systems.
In the present paper we use different approach. We work with elementary theory of continued fractions and apply A. Weil's bounds for Kloostermann
sums.
 Our continued fractions' consideration is connected with a paper by
M. Laurent and A. Nogueira \cite{loor}.

{\bf 3. Lemmata.}

We consider Euclidean plane $\mathbb{R}^2$ with coordinates $(u,v)$. 
Suppose that $|\alpha\delta - \beta \gamma | = 1 $ and put
$$
M = 4 \max ( |\alpha|, |\beta|, |\gamma|, |\delta|, |\delta|^{-1}) 
$$
Consider the lattice
$$
\Lambda =
\left(
\begin{array}{cc}
\alpha &\gamma\cr
\beta&\delta
\end{array}\right)
\mathbb{Z}^2.
$$
The following lemma is a simple result from continued fractions'  theory.

{\bf Lemma 1.}\,\, {\it Suppose that $\beta/\delta \not\in \mathbb{Q}$.
Let $ p_\nu/q_\nu $ and $p_{\nu+1}/ q_{\nu+1}$ be two consequtive convergent fractions to $\beta/\delta$.
Then

{\rm (i)}  
$\left|
\begin{array}{cc}
q_\nu &q_{\nu+1}\cr
p_{\nu}&p_{\nu+1}
\end{array}\right| = (-1)^{\nu}$;

{\rm (ii)} vectors
$$
{\bf e}_1 =    \left(
\begin{array}{c}
\alpha\cr
\beta
\end{array}\right)q_\nu-
\left(\begin{array}{c}
\gamma\cr
\delta
\end{array}\right) p_\nu,\,\,\,\,
{\bf e}_2 =    \left(
\begin{array}{c}
\alpha\cr
\beta
\end{array}\right)q_{\nu+1}- 
\left(\begin{array}{c}
\gamma\cr
\delta
\end{array}\right) p_{\nu+1}
$$
form a basis of $\Lambda$;

{\rm (iii)} 
the set 
$$
\Pi_\nu  =
\{ (u,v) \in \mathbb{R}^2:\,\,\, |u| \le M q_{\nu+1},\,\,\,
|v| \le Mq_{\nu+1}^{-1}\} \subset \mathbb{R}^2
$$
contains a fundamental domain with respect to $\Lambda$,
and hence any its shift $ {\bf e} + \Pi_\nu,\,\, {\bf e} \in \mathbb{R}^2$ contains at least one point of $\Lambda$.}

{\bf Corollary 1.}\,\,{\it
For any  $ \eta$ and $R \in \mathbb{Z}_+$ the set
$$
\Pi_\nu(\eta; R)=
\{ 
 (u,v) \in \mathbb{R}^2:\,\,\, |u | \le 2M R q_{\nu+1},\,\,\,
|v-\eta| \le 2MRq_{\nu+1}^{-1}\}$$
contains $ (2R+1)^2$ points of the lattice $\Lambda$ of the form
\begin{equation}\label{lat}
{\bf e}_0 + l_1{\bf e}_1+l_2{\bf e}_2,\,\,\,\, l_j \in \mathbb{Z},\,\,\, |l_j |\le R,
\end{equation}
with some ${\bf e}_0 \in \Lambda$.

}

{\bf Lemma 2.}\,\,{\it
 There exist $  {\bf  f}_1, {\bf f}_2 \in \Lambda$ such that

{\rm (i)}  $  {\bf  f}_1, {\bf f}_2$ form a basis of $\Lambda$;

{\rm (ii)} $  {\bf  f}_j \in \Pi_\nu ( \eta_j, q_{\nu+1}^r),\,\, j = 1,2$.}

Proof. Consider lattice points of the form (\ref{lat}) from Corollary 1 for parameters $\eta_1$ and $\eta_2$. 
Let these lattice points be
\begin{equation}\label{i1}
{\bf e}_0^1 + l_1^1{\bf e}_1+l_2^1{\bf e}_2
\end{equation}
for parameter 
 $\eta_1$ and
\begin{equation}\label{i2}
{\bf e}_0^2 + l_2^1{\bf e}_1+l_2^2{\bf e}_2
\end{equation}
for parameter 
 $\eta_2$.
We may suppose that $\eta_1 \neq 0$.
Put $ \rho = |\eta_2/\eta_1|$, $  Q =q_{\nu+1}, R = Q^r = q_{\nu+1}^r$ and
$$\left(
\begin{array}{cc}
A_1,&A_2\cr
B_1&B_2
\end{array}\right)
=
\left(
\begin{array}{cc}
q_\nu &q_{\nu+1}\cr
-p_{\nu}&-p_{\nu+1}
\end{array}\right)^{-1}
 \left(\begin{array}{cc}
0&0\cr
\eta_1&\eta_2
\end{array}\right)
.
$$
Then (\ref{rho}) is  valid and $A_1 \asymp Q$.
Under the multiplication by the matrix $\left(
\begin{array}{cc}
q_\nu &q_{\nu+1}\cr
-p_{\nu}&-p_{\nu+1}
\end{array}\right)^{-1}$
lattice points (\ref{i1}) and (\ref{i2}) turn into integer points
$\left(
\begin{array}{c}
x\cr
z
\end{array}\right)$ and
$\left(
\begin{array}{c}
y\cr
w
\end{array}\right)$
respectively, satisfying (\ref{ii}).
By the definition of $r$  there exist $x,y,z,w$ satisfying (\ref{dee}).
So there exist $ l_i^j, \, i,j = 1,2$ such that two points (\ref{i1}) and (\ref{i2}) form a basis of $\Lambda$.$\Box$

{\bf Lemma 3.}\,\,{\it Suppose that
\begin{equation}\label{1}
\left|
\begin{array}{cc}
\xi_1 &\xi_2\cr
\eta_1&\eta_2
\end{array}\right|
=
\left|
\begin{array}{cc}
\Xi_1 &\Xi_2\cr
H_1&H_2
\end{array}\right|
=1.
\end{equation}
Suppose that
\begin{equation}\label{2}
\max_{j=1,2}|H_j -\eta_j |\le \varepsilon
\end{equation}
and $H_1,H_2 \neq 0$. Then
\begin{equation}\label{3}
\left|
\frac{\Xi_1-\xi_1}{H_1}-
\frac{\Xi_2-\xi_2}{H_2}\right|\le  \frac{ \varepsilon \cdot  (|\xi_1|+|\xi_2|)}{|H_1H_2|}.
\end{equation}
}

Proof.
As
$$
\xi_1 \eta_2 - \Xi_1 H_2 - (\xi_2 \eta_1 -\Xi_2 H_1) =0,
$$
from  (\ref{1}) we seee that
$$
|(\xi_1- \Xi_1) H_2 - (\xi_2  -\Xi_2) H_1|\le  \varepsilon \cdot  (|\xi_1|+|\xi_2|).
$$
So 
(\ref{3}) follows.$\Box$

{\bf 4. Proof of theorems.}

Proof of Theorem 1.

 We take large $\nu$.  
Let vectors ${\bf f}_j$ from Lemma 2 be of the form
${\bf f}_j =  \left(\begin{array}{c}\Xi_j\cr H_j\end{array}\right), j = 1,2.$
Then $H_1, H_2 \neq 0$ and
$
\max_{j=1,2}|H_j -\eta_j |\le 2M q_{\nu+1}^{r-1}
$
and
$
\max_{j=1,2}|\Xi_j |\le 2Mq_{\nu+1}^{1+r}
.
$
Put
$
t_\nu = 
\frac{\xi_1-\Xi_1}{H_1} 
.$
Then  by (\ref{3}) we have
$$
\left|t_\nu - \frac{\xi_2 -\Xi_2}{H_2}\right| = O(q_{\nu+1}^{r-1})
.
$$
Of course   we have $ |t_\nu|  = 0(q_{\nu+1}^{1+r})$.
Now
$$
 \left(
\begin{array}{cc}
1 &t_\nu\cr
0&1
\end{array}\right)
\left(
\begin{array}{cc}
\alpha &\gamma\cr
\beta&\delta
\end{array}\right)
\left(
\begin{array}{cc}
l_1^1 &l_1^2\cr
l_2^1&l_2^2
\end{array}\right)=
\left(
\begin{array}{cc}
1 &t_\nu\cr
0&1
\end{array}\right)
\left(
\begin{array}{cc}
\Xi_1 &\Xi_2\cr
H_1&H_2
\end{array}\right)
,
$$
and
$$
 \left(
\begin{array}{cc}
1 &t_\nu\cr
0&1
\end{array}\right)
\left(
\begin{array}{cc}
\Xi_1 &\Xi_2\cr
H_1&H_2
\end{array}\right)-
\left(
\begin{array}{cc}
\xi_1 &\xi_2\cr
\eta_1&\eta_2
\end{array}\right)
= O(q_{\nu+1}^{r-1})
$$
So Theorem 1 follows.$\Box$

Proof of Theorem 2.

Given real $U \ge 1$
put $U_* = \rho (U)$.
By Minkowski's convex body theorem  and condition (\ref{zabo})
one can take primitive point $(q,p)\in \mathbb{Z}^2$ such that
\begin{equation}\label{dia}
U_*\le  q \le U.\,\,\,\,\,|(\beta /\delta ) \cdot q  - p|\le U^{-1}.
\end{equation}
Then this point may be completed to a basis   of $\mathbb{Z}^2$ by a point 
$(q',p')\in \mathbb{Z}^2$ such that
\begin{equation}\label{dia1}
U_*\le  q' \le 2U.\,\,\,\,\,|(\beta/\delta )\cdot q -  p|\le U_*^{-1}.
\end{equation}
From (\ref{dia},\ref{dia1}) we see that the rectangle
$$
 \{ (u,v) \in \mathbb{R}^2:\,\,\, |u| \le MU ,\,\,\,
|v| \le MU_*^{-1}\} \subset \mathbb{R}^2
$$
contains a fundamental domain for $\Lambda$.
Now we follow the argument of the proof of Theorem 1.
We see that $t$ may be taken to be $\ll U^{1+r}$, and we establish the bound
$$
\left|\left|
\left(
\begin{array}{cc}
1 &t\cr
0&1
\end{array}\right)
\left(
\begin{array}{cc}
\alpha &\gamma\cr
\beta&\delta
\end{array}\right)
\left(
\begin{array}{cc}
l_1^1 &l_1^2\cr
l_2^1&l_2^2
\end{array}\right)
-
\left(
\begin{array}{cc}
\xi_1 &\xi_2\cr
\eta_1&\eta_2
\end{array}\right)
\right|\right|\ll U^{r}U_*^{-1}.
$$
By puting $U =T^{\frac{1}{1+r}}$ we have
$$U^{r}U_*^{-1} = T^{\frac{r}{1+r} }/ \rho \left(T^{\frac{1}{1+r}}\right),
$$
and Theorem 2 follows.$\Box$

{\bf 5. About admissible value of $r$.}

Here we show that
any value $r> 3/4$ is good for our purpose.
First of all  we may suppose that $ w=p$ is a prime number
(here we use a well-known  fact that between $Q $ and $Q+Q^{3/4}$ for large $Q$
there exists a prime number,  see \cite{Praa}).

Then we apply the well-known fact that for any two intervals $I_1,I_2$ of lengths $\gg p^{3/4+\varepsilon}$ there exist
$ x\in I_1, y\in I_2$ such that  $xy\equiv 1\pmod{p}$ (see \cite{spa,ou}).
This result follows from A. Weil's bounds for Kloostermann sums.
Now $z  =\frac{xy-1}{p}$ will be an integer.
Easy calculation shows that $x,y,z,w$ will satisfy (\ref{ii}).


\begin{thebibliography}{1}
 

\bibitem{cas} 
J. W. S. Cassels,\,\,
An introduction to the Geometry of Numbers, Springer-Verlag, 1959.

\bibitem{loor}
M. Laurent, A. Nogueira,\,\,
Approximation to points in the plane by  $\rm{ SL}_2(\mathbb{Z})$-orbits,
 Journal of the London Mathematical Society, V. 85, No.2 (2012) pp. 409 - 429.

\bibitem{MW}
 F. Maucourant,
B. Weiss,
\,\,
Lattice actions on the plane revisited,
Geom Dedicata (2012) 157, 1 - 21.

\bibitem{baak} D.W. Morris,\,\,
Ratner's theorem on Unipotent flows,
Chicago lectures in Mathematics Series, 2005.


 \bibitem{Praa}
K. Prachar,\,\,
Primzahlenvertung, Springer-Verlag, 1957.

\bibitem{R1} M. Ratner,\,\,
Raghunathan's topological conjecture and distribution of unipotent flows,
Duke Math. J. 63 (1991), no. 1, 235 - 280.
 

 

\bibitem{spa}
I. Shparlinski,\,\, Modular Hyperbolas,
preprint available at 
 arXiv:1103.2879v2 (2011).

\bibitem{SSS} A. Str\"ombersson,\,\,
On the deviation of ergodic averages for horocycle 
flows, (2003)
preprint,
available at  http://
www.math.uu.se/~astrombe/papers/iha.pdf



\bibitem{ou}
A.V. Ustinov,\,\,
Applications of Kloosterman Sums in Arithmetic and Geometry. (In Russian). LAMBERT Academic Publishing, 2011. 





\end{thebibliography}
 \end{document}